# The Diophantine Equation

$$x^n + y^m = c \cdot x^k \cdot y^\ell,$$

n, m, k, $\ell$, c natural numbers

*Konstantine Zelator*
*Department of Mathematics*
*301 Thackery Hall*
*University of Pittsburgh,*
*Pittsburgh, PA 15260, USA*

*Also: P.O. Box 4280*
*Pittsburgh, PA 15203, USA*

email address:   1) konstantine_zelator@yahoo.com
                 2) kzet159@pitt.edu

## 1. Introduction

The subject matter of this work is the two-variable Diophantine equation,

$$\begin{cases} x^n + y^m = c \cdot x^k \cdot y^\ell \\ \text{where } n, m, k, \ell \text{ and } c \text{ are} \\ \text{natural or positive intergers} \end{cases} \quad \textbf{(1)}$$

We will investigate this Diophantine equation, from the point of view of positive integer solutions. So, for the sake of this article, when we use the term "solution"; we will automatically refer to an ordered pair $(a,b)$ in $\mathbb{Z}^+ \times \mathbb{Z}^+$, which is a solution to **(1)**; that is, when $x = a$ and $y = b$ are substituted into the above equation, then **(1)** becomes a true statement. A preliminary examination of sources such as L.E. Dickson's monumental book, *"History of the Theory of Numbers, Vol. II"* (see [1]), reveals that little or no material can be found pertaining to the Diophantine equation **(1)**. There are a few Diophantine equations that bear some resemblance to **(1)** (for example such is the $x^m X^n + y^m Y^n = z^m Z^n$ equation studied by V. Bouniakowsky, see [1]). The same can be said about reference [2], which is treasure source of number theory material; but it contains no material on equation **(1)**.

A few words about the organization of this paper. In the next section, Section 2, we explain how the material of this work is organized according to eight cases. In Section 3, five results are listed. The first four Results are very well known and are listed without proof. We offer a proof to Result 5 though. As the reader will quickly realize, Result 5 is used repeatedly throughout this work, and it is a central tool in obtaining most of the eight results (Theorems 1 through 8) pertaining to Diophantine equation **(1).** In Sections 5 through 12; the eight cases (described in Section 2) are worked out; one case per section, starting with Section 5 and ending with Section 12. The results obtained for each case, are explicated in the corresponding theorem. So, Theorem 1 contains the results obtained in Case 1; and so on, Theorem 8 contains the results obtained in Case 8. Also, in Section 4 we establish a lemma, Lemma 1.

**2. The Eight Cases**

Looking at equation **(1)**, we see that if d is the greatest common division of *x* and *y*; the **(1)** is equivalent to (by setting $x = dx_1, \ y = dy_1$),

$$\begin{cases} d^n \cdot x_1^n + d^m \cdot y_1^m = c \cdot d^{k+\ell} \cdot x_1^k \cdot y_1^\ell \\ x = d \cdot x_1, \ y = d \cdot y_1 \\ \text{where x, y, d are positive integers} \\ \text{such that } x_1 \text{ and } y_1 \text{ are relatively prime;} \\ \gcd(x_1, y_1) = 1 \end{cases} \quad (2)$$

It will become evident when we start analyzing the various cases which arise from **(2)**; that the key numbers are the integers $n, m,$ and $k + \ell$. First observe that we do not need to distinguish between the grand cases n<m and n>m; for if one switches *x* with *y* in **(1)** (and thus in **(2)**); an equation of the same general form is obtained. Accordingly, there are really eight cases to consider. These are:

Case 1:   $n = m = k + \ell$

Case 2:   $n < m < k + \ell$

Case 3:   $k + \ell < n < m$

Case 4:   $n < k + \ell < m$

Case 5:   $n = k + \ell < m$

Case 6:   $n < m = k + \ell$

Case 7:   $m = n < k + \ell$

Case 8:   $k + \ell < m = n$

## 3. Five Results From Number Theory

The first four results are very will known in elementary number theory.

**Result 1:** Let a and b be relatively prime positive integers, gcd(a,b)=1. And let $n_1$, $n_2$ be natural numbers. Then the integers $a^{n_1}$ and $b^{n_2}$ are relatively prime; that is $\gcd(a^{n_1}, b^{n_2}) = 1$.

**Result 2:** Let a, b, c be positive integers with a and b being relatively prime; gcd(a,b)=1. Then, if a is a divisor of the product b•c; a is a divisor of the integer c.

**Result 3:** Let a and b be positive integers, and let d be their greatest common divisor; $a=da_1$, $b=db_1$; where $a_1$ and $b_1$ are relatively prime positive integers; $\gcd(a_1, b_1)=1$. Also, let $\ell$ be the least common multiple of a and b. Then, (i) $ab = \ell d$ and $\ell = da_1 b_1$. (ii) If d=1 then $a=a_1$, $b=b_1$, and $\ell = ab$.

**Result 4:** Suppose that a and b are relatively prime positive integers; gcd(a,b)=1. Also, assume that $a \mid c$ and $b \mid c$; that is a and b are divisors of some positive integer c. The $ab \mid c$, i.e. the product ab is a divisor of c.

The following result, Result 5, plays a central role in this work. A proof of part *(ii)* of this result can be found in Reference [2].

**Result 5:** Let a and b be positive integer; d the greatest common divisor of a and b; d=gcd(a,b). Also, $a=da_1$, $b=db_1$, where the positive integers $a_1$ and $b_1$ are relatively prime; $\gcd(a_1 b_1)=1$. Consider the two-variable Diophantine equation $x^a = y^b$. Also let $\mathbb{Z}^+$ be the set of positive integers. (i) The solution set S to this Diophantine equation is given by $S = \left\{ (x, y) \mid (x, y) = (t^{b_1}, t^{a_1}); t \in \mathbb{Z}^+ \right\}$. Recall from Result 3 that $b_1 = \dfrac{\ell}{a}$ and $a_1 = \dfrac{\ell}{b}$, where $\ell$ is the least common multiple of a and b. (ii) If a and b are relatively prime; i.e, d=1; then the solution set is $S = \left\{ (x, y) \mid (x, y) = (t^b, t^a); t \in \mathbb{Z}^+ \right\}$.

**Proof**

(i) First, observe that for every positive integer $t$; the pair $(x, y) = (t^{b_1}, t^{a_1})$ is a solution to the said Diophantine equation. Indeed, we have $x^a = (t^{b_1})^a = t^{b_1 a}$; and by Result 3 we know that $\ell$ is the least common multiple $a$ and $b$. Thus $x^a = t^{b_1 a} = t^{b_1 a_1 d} = t^\ell$. Likewise, $y^b = (t^{a_1})^b = t^{a_1 b} = t^{a_1 b_1 d} = t^\ell$. Hence, $x^a = y^b$. Next, we prove the converse, that every solution $(x,y)$ to the equation $x^a = y^b$ must satisfy $(x, y) = (t^{b_1}, t^{a_1})$, for some positive integer $t$. Let $(x,y)$ be a positive integer solution to the Diophantine equation,

$$x^a = y^b \qquad (3)$$

If $x=1$, then $y=1$; and conversely if $y=1$, then $x=1$. So clearly the solution $(x,y)=(1,1)$ is represented in the above form with $t=1$. Assume $x \geq 2$ and $y \geq 2$. By inspection we see that if a prime $p$ divides $x$, then $p$ divides $y$; and conversely, if $p | y$ then $p | x$. This then shows, that $x$ and $y$, have the same prime bases, when they are factored into prime powers (Fundamental Theorem of Arithmetic). So let $p_1, \ldots, p_r$ be the prime bases for $x$ and $y$:

$$\left( x = p_1^{e_1} \ldots p_r^{e_r} \text{ and } y = p_1^{f_1} \ldots p_r^{f_r} \right) \qquad (4)$$

where the exponents $e_1, \ldots, e_r, f_1, \ldots, f_r$ are positive integers. It is clear then from **(4)** and **(3)** that we must have $a \cdot e_i = b \cdot f_i$, for $i = 1, \ldots, r$. And using $a = da_1$, $b = db_1$ we further obtain $\quad a_1 \cdot e_i = b_1 \cdot f_i$, for $i=1,\ldots r$ \qquad (5)

Since $\gcd(a_1, b_1)=1$, by **(5)** and Result 2, it follows that $a_1$ must be a divisor of $f_i$ and $b_1$ a divisor of $e_i$. In effect, we must have, $e_i = k_i \cdot b_1$

$$\begin{pmatrix} e_i = b_1 \cdot k_i \text{ and } f_i = a_1 \cdot k_i, \\ \text{for some positive integers } k_i \\ \text{for } i = 1, \ldots r \end{pmatrix} \qquad (6)$$

Going back to **(4)** we see that $x = \left(p_1^{k_1} \ldots p_r^{k_r}\right)^{b_1}$ and $y = \left(p_1^{k_1} \ldots p_r^{k_r}\right)^{a_1}$; which proves that $(x, y) = (t^{b_1}, t^{a_1})$; $t_1 = p_1^{k_1} \ldots p_r^{k_r}$. This proof is complete. ∎

(ii) This is an immediate consequence of part (i). When d=1, a=$a_1$ and b=$b_1$. ∎

## 4. A Lemma and its Proof

**Lemma 1:** *Consider the Diophantine equation (1) under the condition gcd(x,y)=1. Then, if c=2, equation has the unique solution (x,y)=(1,1). And if c≠2, equation (1) has no solution.*

**Proof**

Equation **(1)** is equivalent to **(2)**. So when condition gcd(x,y)=1=d is imposed, we obtain from **(2)** the equation

$$x_1^n + y_1^m = c \cdot x_1^k \cdot y_1^\ell \quad (7)$$

Clearly, **(7)** implies that $x_1$ divides $y_1^m$; but $x_1$ and $y_1$ are relatively prime; and so by Result 1, $\gcd(x_1, y_1^m) = 1$. But then $(x_1 | y_1^m$ and $\gcd(x_1, y_1^m) = 1)$ clearly imply that $x_1$=1 (since if $x_1$ contained a prime divisor $p$; the $p | y_1^m$ and so $p | y_1$; violating the condition gcd($x_1,y_1$)=1). Likewise, equation **(7)** shows that $y_1 | x_1^n$; and since $\gcd(y_1, x_1^n) = 1$, we conclude that $y_1$=1. Thus $x_1$=$y_1$=1 and by **(7)** we get c=2. Therefore when c=2, equation **(7)** has the unique solution ($x_1,y_1$)=(1,1). And since d=1, x=d$x_1$=1•1=1 and y=d•$y_1$=1•1=1. Equation **(1)**, under the condition gcd(x,y)=1, had the unique solution (x,y)=(1,1). And when c≠ 2, obviously **(7)** has no solutions; and the equation **(1)** has no solution either.

## 5. Case 1: n=m=k+$\ell$

By **(2)**, we obtain $x_1^n + y_1^m = c \cdot x_1^k \cdot y^\ell$; with $\gcd(x_1, y_1) = 1$ which according to the proof of Lemma 1 has the unique solution ($x_1,y_1$)=(1,1) if c=2; otherwise it has no solutions.

Going back to **(2)** we obtain x=d, y=d in the case of c=2. We have the following theorem.

**Theorem 1:** *Suppose that the positive integers* $n, m, k, \ell$ *satisfy the condition* $n = m = k + \ell$. *Then Diophantine equation* **(1)** *has no solution if c≠2, all the solutions of Diophantine equation* **(1)** *are given by* (x,y)=(d,d); *where d can be any positive integer.*

## 6. Case 2: $n < m < k + \ell$

**Theorem 2:** *Let* $n, m, k,$ *and* $\ell$ *be positive integers satisfying the conditions* $n \leq k$ *and* $n < m < k + \ell$. *Then if c≠2, the Diophantine equation* **(1)** *has no solution in positive integers x and y. If on the other hand c=2 Diophantine equation* **(1)** *has the unique solution* (x,y)=(1,1).

**Proof:** Since equation **(1)** is equivalent to **(2)** and $n < m < k + \ell$; equation **(2)** yields

$$x_1^n + y_1^m \cdot d^{m-n} = c \cdot x_1^k \cdot y^\ell \cdot d^{k+\ell-n} \qquad (8)$$

By inspection we see from equation **(8)** that the integer $y_1$ must divide $x_1^n$; and since $\gcd(x_1 y_1) = \gcd(x_1^n, y_1) = 1$, it follows that $y_1=1$. Consequently **(8)** gives

$$x_1^n + d^{m-n} = c \cdot x_1^k \cdot d^{k+\ell-n} \qquad (8a)$$

Since n≤k; we have min{k,n}=n. By **(8a)** it follows that $x_1^n$ must be a divisor of $d^{n-m}$:

$$x_1^n \mid d^{n-m} \qquad (8b)$$

Moreover, the condition $m < k + \ell$ implies that m-n<k+ℓ-n.

And so $\min\{m - n, k + \ell - n\} = m$. Thus **(8a)** implies that

$$d^{n-m} \mid x_1^n \qquad (8c)$$

Since, by **(8b)** and **(8c)**, the positive integers are divisors of each other; they must be equal: $$x_1^n = d^{n-m} \quad (8d)$$

Let $L$ be the least common multiple of the positive integers $n$ and $m-n$. By Result 5, it follows from **(8d)** that $x_1 = t^{\frac{L}{n}}$ and $d = t^{\frac{L}{m-n}}$, for some positive integer $t$. Using these formulas (in terms of $t$) for $x_1$ and $d$ in equation **(8a)**, we obtain

$$2t^L = c \cdot t^{\left(\frac{L}{n} \cdot k + \frac{L}{m-n} \cdot (k+\ell-n)\right)} \quad (8e)$$

Compare the two exponents of $t$: The exponent of $t$ on the right-hand side of **(8e)** is the larger one: Indeed, since $1 \leq n \leq k$; we have $\frac{L}{n} \cdot k = L \cdot \frac{k}{n} \geq L$. And since $\frac{L}{m-n}$ is a positive integer (since $L$ is divisible by $m-n$; as the common multiple of $n$ and $m-n$) and $k+\ell-n$ is also a positive integer (clearly by virtue of $n < m < k+\ell$); it is clear that $\frac{L}{n} \cdot k + \frac{L}{m-n} \cdot (k+\ell-n) \geq L+1 > L$. Therefore **(8e)** is equivalent to,

$$2 = c \cdot t^{\left(\frac{L}{n} \cdot k + \frac{L}{m-n} \cdot (k+\ell-n) - L\right)} \quad (8f)$$

As we showed above, the exponent $\frac{L}{n} \cdot k + \frac{L}{m-n} \cdot (k+\ell-n) - L$ is a positive integer. If this exponent is greater than or equal to 2, then it is clear from **(8f)** that no prime can divide $t$; so in this case $t$ must equal 1, and so by **(8f)**, $c=2$. One obtains, from $x_1 = t^{\frac{L}{n}}$ and $d = t^{\frac{L}{m-n}}$, $x_1 = d_1 = 1$. And since $y_1 = 1$; we get from **(2)**, $x=y=1$; the unique solution $(x,y)=1$. Now if the exponent of $t$ in **(8f)** is equal to 1, we must have

$$1 = \frac{L}{n} \cdot k + \frac{L}{m-n} \cdot (k+\ell-n) - L \Leftrightarrow n \cdot (m-n) = L\left[(m-n)k + n(k+\ell-n) - (m-n)n\right]$$
$$(L+1)n(m-n) = L \cdot \left[(m-n)k + n(k+\ell-n)\right]$$

$$\Leftrightarrow \quad \frac{L+1}{L} = \frac{k}{n} + \frac{k+\ell-n}{m-n} \quad \textbf{(8g)}$$

If L≥2, then $\frac{L+1}{L} = \frac{1}{L} + 1$ falls strictly between 1 and 2:

$1 < 1 + \frac{1}{L} < 2$. But $\frac{k}{n} \geq 1$ (in view of 1≤n≤ k) and $\frac{k+\ell-n}{m-n} > 1$ (in view of

$n < m < k+\ell$; $m-n < k+\ell-n$). And so, if L≥2, the right-hand side of **(8g)** is greater

than 2; but the left-hand side of **(8g)** is less than 2. Impossible.

Furthermore, if L=1; then since L is the least common multiple of n and m-n; it follows

that n=1 and m-n=1; and so m=2. From **(8g)**, with L=1=n and m=2, we obtain

$2 = k + k + \ell - 1$; $3 = 2k + \ell$ which implies, since both k and $\ell$ are positive integers; that

$k = \ell = 1$. However this contradicts the part of the hypothesis: namely that

$m < k+\ell$, obviously with m=2 and $k = \ell = 1$, this hypothesis cannot be satisfied. The

proof is complete. □

## 7. Case 3: $k + \ell < n < m$

**Theorem 3:** *Suppose that the positive integers* $k, \ell, n, m$ *satisfy the condition*

$k + \ell < n < m$.

*(i) If (x,y), and only if, is a solution to the Diophantine equation (1), then*

$$\begin{cases} x = dx_1, \ y = dy_1, \ d = \gcd(x,y), \ d^{m-(k+\ell)} = rx_1^k, \\ d^{n-(k+\ell)} = s \cdot y_1^\ell, \ c = x_1^{n-k} \cdot r + y_1^{m-\ell} \cdot s, \\ \text{where s and r are positive integers; and } x_1, y_1 \text{ are} \\ \text{relatively prime positive integer.} \end{cases}$$

*Note:* Since the positive integer c *is fixed, it is clear from part (i)*

*that the Diophantine equation* **(1)** *has finitely many solutions.*

*(ii) If* c=1, 3, 4, *or* 5; *Diophantine equation* **(1)** *has no positive integer solutions. If* c=2, *equation* **(1)** *has the unique solution* (x,y)=(1,1).

*(iii) If the natural number* c *is the form,* $c = \rho^{m-(k+\ell)} + \rho^{n-(k+\ell)}$, *where* $\rho$ *is a natural number. Then* (x,y)=($\rho,\rho$) *is a solution to equation* **(1)**.

**Proof**: (ii) We start with equation **(2)**, which is of course equivalent to **(1)**. Since $k + \ell < n < m$; the numbers $n - (k+\ell)$ and m-(k+$\ell$) are positive integers and hence **(2)** is equivalent to

$$x_1^n \cdot d^{n-(k+\ell)} + y_1^m \cdot d^{m-(k+\ell)} = c \cdot x_1^k \cdot y_1^\ell \qquad (9)$$

From $k + \ell < n$, it is also clear that the natural number $k$ is smaller than n:1≤k<n. Hence **(9)** is implies that $x_1^k \mid y_1^m \cdot d^{m-(k+\ell)}$; and since $\gcd(x_1^k, y_1^m) = 1$ (since $x_1$ and $y_1$ are relatively prime and by Result 1), it follows that $x_1^k$ must be a divisor of $d^{m-(k+\ell)}$ (by Result 2):

$$d^{m-(k+\ell)} = r \cdot x_1^k, \qquad (9a)$$
for some positive integer r

Likewise, the inequality of $k + \ell < m$ implies $1 \le \ell < m$, and so a similar argument establishes that $y_1^\ell$ must be a divisor of the product $x_1^n \cdot d^{n-(k+\ell)}$; and since $\gcd(y_1^\ell, x_1^n) = 1$; we infer that

$$d^{n-(k+\ell)} = s \cdot y_1^\ell, \qquad (9b)$$
for some positive integer r

Combining **(9)**, **(9a) and (9b)** leads to the condition

$$c = x_1^{n-k} \cdot r + y_1^{m-\ell} \cdot s \qquad \textbf{(9c)}$$

We have shown that the conditions **(9a), (9b), (9c)** are necessary in order for *(x,y)* to be a solution of Diophantine equation **(1)**. These conditions are also sufficient: if **(9a), (9b), (9c),** are satisfied; one easily sees that equation **(9)** can be derived, and hence equations **(2)** and **(1)** as well. We omit the details.

(ii) The condition $k + \ell < n < m$ implies $\ell < n - k$ and k<m-$\ell$. Thus, since $\ell$ and *k* are natural numbers, it follows that

$$(n - k \geq 2 \text{ and m-}\ell \geq 2) \qquad \textbf{(9d)}$$

If c=1, obviously **(9c)** cannot be satisfied, so there are no solutions. If c=2, then **(9c)** necessitates that $x_1^{n-k} \cdot r = 1 = y_1^{m-\ell} \cdot s$, from which we get $x_1 = r = y_1 = s = 1$. And thus from **(9a) and (9b)** we have d=1. Equation **(1)** has the unique solution (x,y)=(1,1). If c=3, then either $x_1^{n-k} \cdot r = 2$ and $y_1^{m-\ell} \cdot x = 1$, or alternatively, $x_1^{n-k} \cdot r = 1$ and $y_1^{m-\ell} \cdot s = 2$. In the first case, given that n-k$\geq$2; we must have $x_1$=1 and r=2; and $y_1$=s=1. But then **(9b)** implies d=1; contradicting **(9a)** in view of r=2. In the latter case, given that $m - \ell \geq 2$. We must have $y_1$=2, s=1, and $x_1$=r=1. Again, a contradiction arises from **(9a) and (9b).**

If c=4, there are three possibilities:

$x_1^{n-k} \cdot r = 2 = y_1^{m-\ell} \cdot s$

or $x_1^{n-k} \cdot r = 3$ and $y_1^{m-\ell} \cdot s = 1$

or $x_1^{n-k} \cdot r = 1$ and $y_1^{m-\ell} \cdot s = 3$

Again, in view of **(9d)**; in all three possibilities we must have $x_1$=$y_1$=1. In the first possibility, r=s=2. And since $m - (k + \ell) > n - (k + \ell) \geq 1$, it follows that $m - (k + \ell) \geq 2$;

which renders **(9a)** impossible by virtue of r=2 and $x_1$=1. Finally, suppose that c=5.

There are four possibilities:

Possibility 1: $x_1^{n-k} \cdot r = 4$ and $y_1^{m-\ell} \cdot s = 1$

Possibility 2: $x_1^{n-k} \cdot r = 1$ and $y_1^{m-\ell} \cdot s = 4$

Possibility 3: $x_1^{n-k} \cdot r = 3$ and $y_1^{m-\ell} \cdot s = 2$

Possibility 4: $x_1^{n-k} \cdot r = 2$ and $y_1^{m-\ell} \cdot x = 3$

Again, in all four possibilities given that **(9d)** holds true; we must have $x_1=y_1=1$. This then means that either:   r=4 and s=1;

or   r=1 and s=4

or   r=3 and s=2

or   r=2 and s=3

In all fours cases, and easy check shows that equations **(9a) and (9b)** are incompatible; i.e., they are contradictory.

(iii) Just take $x_1=y_1=1$, $r = \rho^{m-(k+\ell)}$ and $s=\rho^{n-(k+\ell)}$ and d=ρ. Then since

$c = \rho^{m-(k+\ell)} + \rho^{n-(k+\ell)}$; all three conditions **(9a), (9b)), and (9c)** are satisfied. Thus, (x,y)=(ρ,ρ) is a solution to Diophantine equation **(1)**. □

## 8. Case 4:   $n < k + \ell < m$

**Theorem 4:**   *Suppose that the natural numbers* $n, k, \ell,$ *and* m *satisfy the conditions* $n < k + \ell < m$ *and* n≤k. *Then the following are necessary and sufficient conditions for the pair* (x,y) *to be a solution to the Diophantine equation* **(1)**:

*There exist positive integers* d, R, S, *and* $x_1$ *such that*

$$\begin{cases} x = d \cdot x_1, \ y=d (\text{so } y_1 = 1), \ d^{m-(k+\ell)} = R \cdot S, \\ x_1^n = d^{k+\ell-n} \cdot S, \quad S \cdot (1+R) = c \cdot x_1^k, \\ \text{and } d^{m-n} = R \cdot x_1^n \end{cases}$$

**Proof:** From **(2)** and since $n < k + \ell < m$, we obtain

$$x_1^n + y_1^m \cdot d^{m-n} = c \cdot d^{k+\ell-n} \cdot x_1^k y_1^\ell \tag{10}$$

From $k + \ell < m$, it is clear that the positive integer $\ell$ is less than $m$: $\ell < m$. Accordingly, **(10)** implies that $y_1^\ell$ a divisor of $x_1^n$; but $\gcd(x_1^n, y_1^\ell) = 1$, in virtue $\gcd(x_1, y_1)=1$ and Result 1. Hence $y_1^\ell | x_1^n$ implies $y_1=1$. And so equation **(10)** yields

$$x_1^n + d^{m-n} = c \cdot d^{k+\ell-n} \cdot x_1^k \tag{10a}$$

And from **(10a)** we get that since $n \leq k$, $x_1^n | d^{m-n}$; and so we have

$$\begin{cases} d^{m-n} = R \cdot x_1^n \\ \text{for some positive integer R} \end{cases} \tag{10b}$$

Moreover, the hypothesis $k + \ell < m$ implies $k + \ell - n < m - n$. Clearly then equation **(10a)** implies that the integer $d^{k+\ell-n}$ must be a divisor of $x_1^n$:

$$\begin{cases} x_1^n = S \cdot d^{k+\ell-n} \\ \text{for some positive integer S} \end{cases} \tag{10c}$$

From **(10c) and (10b)** it readily follows that

$$d^{m-(k+\ell)} = R \cdot S \tag{10d}$$

Combining **(10a), (10b), and (10c)** produces

$$d^{k+\ell-n} \cdot S + R \cdot S \cdot d^{k+\ell-n} = c \cdot d^{k+\ell-n} \cdot x_1^k$$
$$\text{whence } S \cdot (1+R) = c \cdot x_1^k \tag{10e}$$

13And of course (by **(2)**) $x = x_1 \cdot d$ and $y = y_1 \cdot d = d$ (since $y_1 = 1$). Conditions **(10b), (10c), (10d), (10e),** establish the necessity part of the theorem. To show that they are sufficient, is pretty straightforward. We omit the details. □

The following theorem, is a corollary of Theorem 4:

**Theorem 4a:** *Suppose that* n=2, m=6, *and* k=2=$\ell$; *or alternatively* n=2, m=6, k=3, *and* $\ell$=1. *Then if* c≠2, *the Diophantine equation* **(1)** *has no positive integer solution. And if* c=2, (x,y)=(1,1) *is the unique solution of equation* **(1).**

**Proof:** We apply Theorem 4. First observe that in either case of the hypothesis, $k + \ell = 4$. And so $m - (k + \ell) = 6 - 4 = 2$, n=2, m-n=6-2=4, and so $k + \ell - n = 4 - 2 = 2$. And with *k* being either 2 or 3. And so, n=2≤k≤3 and $n < k + \ell < m$ are satisfied. According to Theorem 4; the pair (x,y) is a solution to **(1)** precisely when,

$$\begin{cases} x = d \cdot x_1,\ y=d,\ d^2 = R \cdot S, \\ x_1^2 = d^2 \cdot S, \quad S \cdot (1+R) = c \cdot x_1^k, \\ \text{and } d^4 = R \cdot x_1^2 \\ \text{where } d,\ x_1,\ R \text{ and } S \text{ are positive integers} \end{cases} \quad \textbf{(10i)}$$

Now, the equation $x_1^2 = d^2 \cdot S$, implies that S must be a perfect or integer square (this follows easily, for example, from the Fundamental Theorem of Arithmetic). Hence,

$$\begin{cases} S = S_1^2 \\ \text{for some natural number } S_1 \end{cases} \quad \textbf{(10ii)}$$

Going back to $d^2 = R \cdot S$, we now have $d^2 = R \cdot S_1^2$, which implies (same argument as before) that

$$\begin{cases} R = R_1^2 \\ \text{for some natural number } R_1 \end{cases} \quad \textbf{(10iii)}$$

Going back to **(10i)** once more we see that with $R = R_1^2$ and $S = S_1^2$ we obtain,

$$\left(d = R_1 \cdot S_1 \text{ and } x_1 = d \cdot S_1 = R_1 \cdot S_1^2\right) \qquad \textbf{(10iv)}$$

From **(10iv)** and $S \cdot (1+R) = c \cdot x_1^k$ in **(10i)** we obtain

$$S_1^2 \cdot (1 + R_1^2) = c \cdot R_1^k \cdot S_1^{2k} \qquad \textbf{(10v)}$$

If both $S_1$ and $R_1$ are greater than 1; since $k \geq 2$, it is clear that $2k-2 \geq 2$; and so, since $c \geq 1$; we obtain $c \cdot S_1^{2k-2} \geq c \cdot 2^{2k-2} \geq c \cdot 2^2 \geq 4$. But then, $4R_1^k > 1 + R_1^2$ (since $k \geq 2$ and $R_1 \geq 2$), rendering **(10v)** impossible. Suppose that $S_1=1$ and $R_1 \geq 2$ in **(10v)**. Then $1 + R_1^2 = c \cdot R_1^k$, with $R_1 \geq 2$. And so if $k=3$, then it is clear that $c \cdot R_1^k \geq R_1^3 > 1 + R_1^2$, since $R_1 \geq 2$. Now, if $k=2$ and $R_1 \geq 2$, we have $1 + R_1^2 = c \cdot R_1^2$, which is not possible when $c=1$. But is also impossible $c \geq 2$, since $c \cdot R_1^2 \geq 2 \cdot R_1^2 > R_1^2 + 1$, which is clearly the case on account of $R_1 \geq 2$. Next assume that $R_1=1$ and $S_1 \geq 2$ in **(10v)**. We have $2 \cdot S_1^2 = c \cdot S_1^{2k}$; $2 = c \cdot S_1^{2k-2}$ which is impossible since $c \cdot S_1^{2k-2} \geq S_1^{2k-2} \geq 2^{2k-2} \geq 2^2 = 4 > 2$. The only remaining possibility is $S_1=1=R_1$. Then equation **(10v)** is satisfied only when $c=2$. And from $S_1=1=R_1$ and equations **(10i)** - **(10iv)** we get $x_1=1$, $d=1$, $R=1$, $S=1$. And so, $(x,y)=(1,1)$ is the unique solution of equation **(1)** in the case of $c=2$. □

## 9. Case 5: $n = k + \ell < m$

**Theorem 5:** *Let m, n, k, $\ell$ be natural numbers such that $n = k + \ell < m$. Consider the Diophantine equation* **(1)**.

*(i) If $c=1$, the Diophantine equation* **(1)** *has no positive integer solutions.*

*(ii) If $c \geq 2$, let M be the maximum or largest positive u integer such that $u^{n-k} < c$ (since $c \geq 2$, such an integer M exists, see remarks below). Then, (x,y) is a solution to the*

*Diophantine equation* **(1)** *if and only if there exist positive integers* $x_1$, r *and* d *such*

*that* $\begin{cases} 1 \leq x_1 \leq M, \ x_1^{n-k} + r = c, \ d^{m-n} = r \cdot x_1^k \\ x = d \cdot x_1, \ \text{and y=d (so that } y_1 = 1) \end{cases}$

*(iii) If* c=2, (x,y)=(1,1) *is the unique solution to equation* **(1)**.

*(iv) If* c=3, n-k$\geq$2 *(note that* n>k *since* n=k+$\ell$ *), and* m-n=1, *then* (x,y)=(2,2) *is the unique solution to equation* **(1)**.

*(v) If* c=3, n-k$\geq$2, *and* m-n$\geq$2, *then equation* **(1)** *has no positive integer solutions.*

*(vi) If* c=3, n-k=1 *(i.e.* $\ell$=1*), and* m-n=1. *Then equation* **(1)** *has exactly two solutions. The pairs* (x,y)=(2,2), ($2^{k+1}, 2^k$).

*(vii) If* c=3, n-k=1 *(i.e.* $\ell$=1*), and* m-n$\geq$2. *The Diophantine equation* **(1)** *has no solutions in positive integers* x *and* y; *unless the integers* n, k, $\ell$, *and* m *satisfy the following three conditions:* k=$\rho \cdot$v, m=$\ell$+v$\cdot\rho$+v, n=$\rho \cdot$v+$\ell$, *for some positive integers* $\rho$ *and* v. *If the integers* n, k, $\ell$, m *satisfy these three conditions then equation* **(1)** *has the unique solution* (x,y)=($2^{\rho+1}, 2^\rho$). *(Note that if* n, k, $\ell$, m *satisfy these three conditions; then* n=k+$\ell$<m, *as required by the overriding hypothesis of the theorem).*

**Remarks:**

1. If c $\geq$2, the integer M exists, since the set

$S = \{u | u \text{ is a positive integer such that } u^{n-k} < c\}$ is nonempty: it contains the number 1: $1^{n-k}$=1<2$\leq$c. The set S is also bounded above since obviously, *u* cannot exceed the value c-1.

2. It is evident from the conditions in part (ii), that when c$\geq$2, Diophantine equation **(1)** has most M solutions.

**Proof**

From equation **(2)** we obtain

$$x_1^n + y_1^m d^{m-n} = c \cdot x_1^k y_1^\ell$$

$$\tag{11}$$

Since $n = k + \ell < m$, we have $\ell < m$. And so **(11)** implies that $y_1^\ell$ is a divisor of $x_1^n$; and since $\gcd(y_1^\ell, x_1^n) = 1$; we infer that $y_1 = 1$. Thus **(11)** yields,

$$x_1^n + d^{m-n} = c \cdot x_1^k \tag{11a}$$

(i) If c=1, then since k<n (by virtue of n=k+$\ell$ ), it is clear that equation **(11a)** is impossible in positive integers $x_1$ and d.

(ii) Suppose c≥2. On account of k<n, we see from **(11a)** that $x_1^k \mid d^{m-n}$ ; therefore,

$$\left\{ \begin{array}{l} d^{m-n} = r \cdot x_1^k \\ \text{for some positive integer } r \end{array} \right\} \tag{11b}$$

Combining **(11a) and (11b)** gives

$$x_1^{n-k} + r = c \tag{11c}$$

And of course we have $1 \leq x_1 \leq M$, where M is the largest positive integer *u* such that $u^{n-k} < c$. Also, it is clear that if **(11c), (11b) and (11a)** are satisfied; then $(x,y) = (x_1 \cdot d, d)$ is a solution to **(1)**.

(iii) If c=2, then **(11c)** necessitates that $x_1 = r = 1$; and so d=1, from **(11b)**. We obtain the unique solution (x,y)=(1,1).

(iv) Suppose that c=3, n-k≥2, and m-n=1. Since n-k≥2 and c=3; **(11c)** gives $x_1$=1 and r=2. And since m-n=1, equation **(11b)** produces d=2. We have x=d·$x_1$=2·1=2 and y=d·$y_1$=2·1=2. the pair (x.y)=(2,2) is the unique solution of equation **(1)**.

(v) Suppose that c=3, n-k≥2, and m-n≥2. As in the previous case, since n-k≥2, equation **(11c)** produce $x_1$=1 and r=2. But **(11b)** gives $d^{m-n}$=2, which is impossible since m-n≥2. No solutions.

(vi) Suppose c=3, n-k=1 (and so $\ell$=1), and m-n=1. From equation **(11c)** we obtain $x_1+r=c=3$; which means either $x_1=2$ and r=1; or $x_1=1$ and r=2. If $x_1=1$ and r=2, the equation **(11b)** becomes $d=r \cdot x_1=r \cdot 1=2$. We obtain the solution (x,y)=(2,2). If on the other hand, $x_1=2$ and r=1, then **(11b)** becomes $d=2^k$. We obtain the solution;

$$x = x_1 \cdot d = 2 \times 2^k = 2^{k+1} \text{ and } y = d = 2^k. \text{ So, } (x,y) = (2^{k+1}, 2^k).$$ Equation **(1)** has two solutions: (x,y)=(2,2), $(2^{k+1}, 2^k)$.

(vii) Finally, assume that c=3, n-k=1 (and so $\ell$=1), and m-n$\geq$2. From equation **(11c)** we know that either $x_1=2$ and r=1; or alternately $x_1=1$ and r=2. In the latter case, equation **(11b)** gives $d^{m-n}=2$, which is impossible with m-n$\geq$2. In the former case, we obtain $d^{m-n}=2^k$, which implies that *d* must be a power of 2: $d=2^\rho$, for some positive integer $\rho$.

Thus $\qquad 2^{\rho(m-n)} = 2^k;$ and hence $\rho$(m-n)=k $\qquad$ **(11d)**

Applying the hypothesis that n=k+$\ell$; we get from **(11d)**

$$\rho(m - \ell) = k \cdot (\rho + 1) \qquad \textbf{(11e)}$$

But $\rho$ and $\rho$+1 are relatively prime: $\gcd(\rho, \rho+1) = 1$. Thus, by Result 2, $\rho$ must be a divisor of *k*:

$$\begin{cases} k = \rho \cdot v, \\ \text{for some positive integer } \rho \end{cases} \qquad \textbf{(11f)}$$

Combining **(11f), (11e), (11d)** and the hypothesis $n = k + \ell$; we obtain the three sought after conditions, namely,

$$(k = \rho \cdot v, \ m = \ell + v \cdot \rho + v, \ n = \rho \cdot v + \ell) \qquad \textbf{(11g)}$$

Recall, from above that $d=2^\rho$, and $x_1=2$, r=1. Hence, when the three conditions in **(11g)** are satisfied, equation **(1)** has the unique solution $x=d \cdot x_1=2^{\rho+1}$, $y=d \cdot y_1=d \cdot 1=d=2^\rho$; (x,y)=$(2^{\rho+1}, 2^\rho)$. Otherwise, equation **(1)** has no solutions. End of proof. $\square$

**10: Case 6:** $n = k + \ell < m$

*Suppose that the positive integers n, k, $\ell$, m satisfy the conditions*

$k + \ell = n > m$ *and* $m \leq \ell$. *Let L be the least common multiple of the positive integers m-n and m; L=l.c.m. (m-n,m). And let* $e_1, e_2, e_3,$ *be the positive integers defined as*

*follows:* $e_1 = \dfrac{L}{m} \cdot (m-n), \quad e_2 = L, \quad e_3 = L + \dfrac{L}{m \cdot \ell}$ *obviously,* $e_2 < e_3$. *Consider the Diophantine equation* **(1)**.

*(i) If $e_1 < e_2 < e_3$, then if c=2, equation* **(1)** *has the unique solution (x,y)=(1,1). Otherwise, for $c \neq 2$,* **(1)** *has no solutions.*

*(ii) If $e_2 < e_1 < e_3$, then for c=2, equation* **(1)** *has the unique solution (x,y)=(1,1). Otherwise, $c \neq 2$,* **(1)** *has no solution in positive integers.*

*(iii) If $e_2 < e_3 \leq e_1$, then for c=2, equation* **(1)** *had the unique solution (x,y)=(1,1); no solutions if $c \neq 2$.*

*(iv) If $e_1 = e_2 = <e_3$, then (as easily seen) n=2m; n-m=m; L=m; and so $e_1=e_2=m$.*

*If c=2, equation* **(1)** *has the unique solution (x,y)=(1,1).*

*If c=1, and $e_3-e_1=1$; then $e_3=m+1$, and so $\ell=1$; k=2m-1. Equation* **(1)** *had the unique solution (x,y)=(2,4).*

*If c=1 and $e_3-e_1 \geq 2$, equation* **(1)** *has no solutions.*

*Likewise, if $c \geq 3$, equation* **(1)** *has no solution in positive integers.*

**Proof:** From **(2)** and in view of $m < n = k + \ell$, we obtain

$$x_1^n \cdot d^{n-m} + y_1^m = c \cdot d^{n-m} \cdot x_1^k \cdot y_1^\ell \qquad (12)$$

Since $n = k + \ell$, the positive integer $k$ is smaller than $n$. Equation **(12)** implies $x_1^k \mid y_1^m$; and since $\gcd(x_1^k, y_1^m) = 1$; it follows that $x_1 = 1$. Thus **(12)** takes the form

$$d^{n-m} + y_1^m = c \cdot d^{n-m} \cdot y_1^{\ell} \qquad (12a)$$

By hypothesis, $m \leq \ell$. Thus **(12a)** shows that $y_1^m$ is a divisor of $d^{n-m}$:

$$y_1^m \mid d^{n-m} \qquad (12b)$$

Likewise from **(12a)**, we see that,

$$d^{n-m} \mid y_1^m \qquad (12c)$$

Hence, by **(12b)** and **(12c)**, it is clear that the positive integers $d^{n-m}$ and $y_1^m$ must be

equal; $\qquad y_1^m = d^{n-m} \qquad (12d)$

By **(12d)** and Result 5 it follows that,

$$\begin{cases} y_1 = t^{\frac{L}{m}} \text{ and } d = t^{\frac{L}{n-m}}; \text{ where } L = \text{l.c.m.}(m, n\text{-}m), \\ \text{and t a positive integer} \end{cases} \qquad (12e)$$

And $e_1 = \dfrac{L}{m} \cdot (n-m)$, $e_2 = L$, $e_3 = L + \dfrac{L}{m} \cdot \ell$, and using equation **(12a)** and also **(12e)**

we get, $\qquad t^{e_1} + t^{e_2} = c \cdot t^{e_3} \qquad (12f)$

(i) If $e_1 < e_2 < e_3$, then clearly equation **(12f)** shows that $t^{e_1} \mid 1$; $t=1$, and so clearly equation

**(12f)** can be satisfied only when $t=1$ and $c=2$; in which case $d=1$, $y_1=1$; and so

$x = x_1 \cdot d = 1 \cdot 1 = 1$, $y = y_1 \cdot d = 1 \cdot 1 = 1$; so we obtain the unique solution $(x,y)=(1,1)$. Otherwise if

$c \neq 2$, there is no solution.

(ii) If $e_2 < e_1 < e_3$; then **(12f)** implies $t^{e_2} \mid 1$; $t=1$. Again, as in (i), we must have $t=1$ and $c=2$;

$(x,y)=(1,1)$. Otherwise, there are no solutions.

(iii) If $e_2 < e_3 \leq e_1$; then $t^{e_2} | t$; t-1, same as in the previous case (ii). Unique solution (x,y)=(1,1) if c=2. No solutions otherwise.

(iv) If $e_1 = e_2 < e_3$; then n-m=m; n=2m, and so L=m and $e_1 = e_2 = m$. If c=2, as usual, (x,y)=(1,1) is the unique solution; since **(12f)** implies $2 = c \cdot t^{e_3 - e_1}$; and so with c=2; we must have t=1. If c=1 and $e_3 - e_1 = 1$. Then $e_3 = 1 + e_1 = 1 + m$; and so from $e_3 = L + \frac{L}{m} \cdot \ell$, we get $e_3 = m + \frac{m}{m} \cdot \ell$; $m + 1 = m + \ell$; $\ell = 1$. And so k=n-$\ell$=2m-1. Also from equation **(12f)** we have, $2 = t^{e_3 - e_1} = t^{m+1-m} = t$; t=2. And so, $y_1 = t^{\frac{m}{m}} = t = 2$, $d = t^{\frac{m}{m}} = t = 2$ (from **(12d)**).

Hence $x = x_1 \cdot d = 1 \cdot 2 = 2$, $y = y_1 \cdot d = 2 \cdot 2 = 4$. The unique solution is the pair (x,y)=(2,4).

Back to equation **(12f)**. Since $e_1 = e_2 < e_3$, we have $2 = c \cdot t^{e_3 - e_1}$. Clearly if c≥3; or if c=1 and $e_- e_1 \geq 2$; the last equation cannot be satisfied for any positive integer t.

End of proof. □

## 11: Case 7: $m = n < k + \ell$

**Theorem 7:** *Suppose that the natural numbers* m, n, k, $\ell$ *satisfy the condition* $m = n < k + \ell$.

*Then, if* c=2, *the Diophantine equation* **(1)** *has the unique solution* (x,y)=(1,1).

*If* c≥3, *equation* **(1)** *has no solutions in positive integers* x *and* y.

*If* c=1 *and in a addition* $k + \ell = n + 1$, *equation* **(1)** *had the unique solution* (x,y)=(2,2).

*If* c=1 *and* $k + \ell - n \geq 2$, *equation* **(1)** *has no solutions.*

**Proof:** From **(2)**, since $m = n < k + \ell$, we obtain

$$x_1^n + y_1^n = c \cdot d^{k+\ell-n} \cdot x_1^k \cdot y_1^\ell \qquad (13)$$

Equation **(13)** implies that $x_1^{\min\{k,n\}}$ is a divisor of $y_1^n$; and since $\gcd(x_1 y_1)=1$; again we conclude that $x_1=y_1=1$. Therefore **(13)** takes the form

$$2 = c \cdot d^{k+\ell-n} \qquad (13a)$$

If $c \geq 3$, obviously **(13a)** cannot be satisfied. And if $c=1$, **(13a)** can be satisfied only if $d=2$ and $k+\ell-n=1$; $k+\ell=n+1$. In which case $x=x_1 \cdot d=1 \cdot 2=2$ and $y=y_1 \cdot d=1 \cdot 1=2$. So $(x,y)=(2,2)$ is the unique solution. Moreover, if $c=1$ and $k+\ell-n \geq 2$, **(13a)** is impossible. Finally if $c=2$, then **(13a)** is satisfied only for $d=1$ in which case $(x,y)=(1,1)$ is the unique solution (since $x_1=y_1=d=1$). □

## 12. Case 8: $k + \ell < m = n$

**Theorem 8:** *Suppose that the positive integer* k, $\ell$, m, n *satisfy the condition* $k + \ell < m = n$. *Consider the Diophantine equation* **(1)**.

*(i) The pair (x,y) is solution the Diophantine equation **(1)** if and only if the following necessary and sufficient conditions are satisfied:*

$$\begin{cases} d = \gcd(x,y), \quad x = d \cdot x_1, y = d \cdot y_1, \quad d^{n-(k+\ell)} = x_1^k \cdot y_1^\ell \cdot r, \\ r \cdot (x_1^n + y_1^n) = c; \quad \text{for positive integers } x_1, y_1, d, \text{ and} \\ r; x_1 \text{ and } y_1 \text{ are relatively prime.} \end{cases}$$

*(ii) If $c=2$, equation **(1)** has the unique solution $(x,y)=(1,1)$.*

*(iii) If $c=1$, equation **(1)** has no positive integer solutions.*

*(iv) If c=4, and* $n = k + \ell$ *equation* **(1)** *had the unique solution* (x,y,)=(2,2). *Otherwise, if* c=4 *and* $n - (k + \ell) \geq 2$; *no solution.*

*(v) If* $c \geq 3$ *and c is a prime and* n *is an odd integer, the Diophantine equation* **(1)** *has no solutions.*

**Proof**:

(i) From equation **(2)** and since $k + \ell < m = n$, we get

$$x_1^n \cdot d^{n-(k+\ell)} t y_1^n \cdot d^{n-(k+\ell)} = c \cdot x_1^k \cdot y_1^\ell \qquad (14)$$

Clearly, by virtue of $k + \ell < n$; we have $1 \leq k < n$. Thus **(14)** shows that $x_1^k \mid y_1^n \cdot d^{n-(k+\ell)}$ ; and since $\gcd(x_1^k, y_1^n) = 1$, it follows that (by Result 2) that

$$x_1^k \mid d^{n-(k+\ell)} \qquad (14a)$$

Likewise, since $\ell < n$ (from $k + \ell < n$) a similar argument establishes that $y_1^\ell$ must be a divisor of $d^{n-(k+\ell)}$:

$$y_1^\ell \mid d^{n-(k+\ell)} \qquad (14b)$$

Since $x_1^k$ and $y_1^\ell$ are relatively prime, it follows by Result 4, **(14a)** and **(14b)** that the product $x_1^k \cdot y_1^\ell$ must be a divisor of $d^{n-(k+\ell)}$:

$$\begin{cases} d^{n-(k+\ell)} = x_1^k \cdot y_1^\ell \cdot r \\ \text{for some positive integer r} \end{cases} \qquad (14c)$$

By equations **(14c)** and **(14)**, it easily follows that

$$r \cdot (x_1^n + y_1^n) = c \qquad (14d)$$

Equations **(14c)** and **(14d)** (together with the usual x=d•x$_1$ and y=d•y$_1$) are the necessary conditions we wished to establish. They are also sufficient: it is an easy matter to show that **(14c)** and **(14d)** imply **(14)**; and hence **(2)**, and equivalently equation **(1)**.

(ii) By inspection, when c=2, **(14d)** gives r=1=x$_1$=y$_1$; and so from **(14c)** we get d=1. Thus x=x$_1$•d=1, y=y$_1$•d=1; (x,y)=(1,1) is the unique solution.

(iii) Clearly when c=1, equation **(14d)** cannot be satisfied. Equation **(1)** has no solution.

(iv) Let c=4. Obviously, it is clear from **(14d)** that r cannot equal 4. If r=1, then from **(14d)** we get $x_1^n + y_1^n = 4$. But n≥3, in virtue of $n > k + \ell \geq 1+1 = 2$. But then with n≥3, clearly the last equation cannot have a solution in positive integers x$_1$ and y$_1$. The only remaining possibility is r=2: we get $x_1^n + y_1^n = 2$ which implies x$_1$=y$_1$=1. Then from **(14c)**, we get $d^{n-(k+\ell)} = 2$; which is possible only when d=2 and n-(k+$\ell$)=1. From d=2, we get x=d•x$_1$=2, y=d•y$_1$=2. Thus (x,y)=(2,2) is the unique solution. Otherwise, if n-(k+$\ell$)≥2, equation **(1)** has no solution.

(v) If c≥3 and c is a prime. Then **(14d)** shows that either r=1 or r=c. Obviously, the possibility r=c is ruled out by inspection. If r=1, the **(14d)** yields

$$x_1^n + y_1^n = c \qquad (14e)$$

However, since n≥3 (by virtue of n>k+$\ell$ )≥1) and n is odd, the left-handed side of **(14e)** factors into two proper divisors of c (i.e. two positive integers strictly between 1 and c):

$(x_1 + y_1)(x_1^{n-1} - x_1^{n-2} \cdot y_1^1 + ... - x_1^1 \cdot y_1^{n-2} + y_1^{n-1}) = c,$ which is impossible since c is a prime and each of the two factors on the left-hand side of the last equation is a positive integer strictly between 1 and c. □